\documentclass[a4paper,10pt]{amsart}
\usepackage{mathrsfs}
\usepackage[arrow,matrix]{xy}
\usepackage{amsmath,amssymb, amscd, amsthm,mathrsfs}
\usepackage{hyperref}
\theoremstyle{plain}
\newtheorem{thm}{Theorem}[section]
\newtheorem{cor}[thm]{Corollary}
\newtheorem{lem}[thm]{Lemma}
\newtheorem{prop}[thm]{Proposition}
\newtheorem{defn}[thm]{Definition}

\newtheorem{rem}[thm]{Remark}

   \def\op{\oplus} \def\ot{\otimes}
\def\Hom{\operatorname {Hom}}\def\D{\mathcal{D}}

\def\Ext{\operatorname {Ext}}

\def\k{\mathbf k}\def\M{\mathcal{M}}

\begin{document}
\title[Calabi-Yau coalgebras]{\bf Calabi-Yau coalgebras}
\author{J.-W. He, B. Torrecillas, F. Van Oystaeyen and Y. Zhang}
\address{J.-W. He\newline \indent Department of Mathematics, Shaoxing College of Arts and Sciences, Shaoxing Zhejiang 312000,
China\newline \indent Department of Mathematics and Computer
Science, University of Antwerp, Middelheimlaan 1, B-2020 Antwerp,
Belgium} \email{jwhe@usx.edu.cn}
\address{B. Torrecillas\newline\indent Department of Algebra and
Analysis, University of Almeria, E-04071, Almeria, Spain}
\email{btorreci@ual.es}
\address{F. Van Oystaeyen\newline\indent Department of Mathematics and Computer
Science, University of Antwerp, Middelheimlaan 1, B-2020 Antwerp,
Belgium} \email{fred.vanoystaeyen@ua.ac.be}
\address{Y. Zhang\newline
\indent Department WNI, University of Hasselt, Universitaire Campus,
3590 Diepenbeek, Belgium} \email{yinhuo.zhang@uhasselt.be}

\date{}

\begin{abstract}  We present a method for constructing the minimal injective
resolution of a simple comodule of a path coalgebra of quivers with
relations. Dual to the Calabi-Yau condition of algebras,  we introduce the
concept of a Calabi-Yau coalgebra. Then we describe the Calabi-Yau coalgebras
of lower global dimensions. An appendix is included for listing
some properties of cohom functors.
\end{abstract}


\maketitle

\section*{Introduction}
Calabi-Yau algebras and categories, because of their links to mathematical physics, algebraical geometry, representation theory, $\dots$ etc., were intensively studied in recent years (cf.\cite{Gin,BT,IR,Boc,CZ}). Bocklandt proved that a graded Calabi-Yau algebra (defined by a finite quiver) of dimension 2 is a preprojective algebra, while a 3-dimensional graded Calabi-Yau algebra is determined by a superpotential (cf.\cite{Boc}). Recently Berger and Taillefer in \cite{BT} defined a special class of nongraded Calabi-Yau algebras, namely the Poincar\'{e}-Birkhoff-Witt deformations of graded Calabi-Yau algebras of dimension 3 . But for general nongraded Calabi-Yau algebras, we don't have much information of them. In fact, it is not easy to study Calabi-Yau the property for general nongraded algebras. Nevertheless, a noetherian complete algebra shares many similar properties with a connected graded algebra (cf. \cite{WZ,CWZ}). Our naive idea is to consider when a noetherian complete algebra is Calabi-Yau. It is well-known that a noetherian complete algebra with cofinite Jacobson
radical is the dual algebra of an artinian coalgebra (cf. \cite{HR}). Since any coalgebra has the locally finite property, we could take this advantage to attack the problem by using coalgebras. However the Calabi-Yau property of coalgebras has not yet been studied. So in this paper we try to lay a foundation for Calabi-Yau coalgebras. In a subsequent paper, we will discuss the Calabi-Yau property of artinian coalgebras (cf. \cite{HTVZ}).

The paper is organised as follows. In Section 1, we introduce the Calabi-Yau condition to coalgebras, and discuss some basic properties. In Section 2, we show that a Calabi-Yau coalgebra of dimensions 0 is exactly a cosemisimple coalgebra. In particular, a semiperfect coalgebra is Calabi-Yau if and only if it is cosemisimple.  A Calabi-Yau coalgebra of dimension 1 is precisely a direct sum of (possibly infinite) copies of the coalgebra $\k[x]$. For Calabi-Yau coalgebras of higher dimensions, we will mainly focus us on the path coalgebras of quivers with
relations. This is reasonable because any coalgebra is
Morita-Takeuchi equivalent to a basic coalgebra (cf. \cite{CG}),
which is a subcoalgebra of the path coalgebra of its Gaberial quiver (cf.
\cite{CM}) if we work over an algebraically closed field. For a path coalgebra of a quiver $(Q,\Omega)$ with relations (cf. \cite{Sim}), we give a construction of the first two steps of the minimal injective resolutions of a simple comodule through the arrows and
the relations of $Q$. This enables us to analyze the necessity
conditions on the quivers and the relations such that the path
coalgebras are Calabi-Yau of dimensions 2 or 3. From there we obtain the results dual to the ones in \cite{Boc}.

Throughout $\k$ is an algebraically closed field of characteristic
zero. All the algebras and coalgebras involved are over $\k$;
unadorned $\ot$ means $\ot_\k$ and Hom means Hom$_\k$.

\section{Calabi-Yau condition for coalgebras}

Let $C$ be a coalgebra. Let $\D^b({}^C\!\M)$ be the bounded derived
category of $C$. Consider the full triangulated subcategory
$\D^b_{fd}({}^C\!\M)$ consisting of complexes with finite
dimensional cohomology. Recall that a $\k$-linear category $\mathcal{T}$ is {\it Hom-finite} if for any two objects $X$ and $Y$ in
$\mathcal{T}$, $\Hom_{\mathcal{T}}(X,Y)$ is a finite
dimensional $\k$-vector space. In general, $\D^b_{fd}({}^C\!\M)$ is not a
Hom-finite $\k$-linear category. Since any simple comodule is finite dimensional, we
have that $\D^b_{fd}({}^C\!\M)$ is Hom-finite if and only if, for
any simple comodules $M$ and $N$, $\Ext^i_C(N,M)$ is finite
dimensional for all $i\ge0$. For example, $\D^b_{fd}({}^C\!\M)$ is
Hom-finite if $C$ is a (left) strictly quasi-finite coalgebra (cf.
\cite{GNT}). Recall that a Hom-finite triangulated category
$\mathcal{T}$ is called a {\it Calabi-Yau category of dimension $n$}
if there are natural isomorphisms $\Hom_{\mathcal{T}}(X,Y)^*\cong
\Hom_{\mathcal{T}}(Y,X[n])$ for all $X,Y\in\mathcal{T}$, that is,
the $n$-th shift functor is a Serre functor (cf. \cite[Appendix]{Boc}).

\begin{defn}\label{cdef1} {\rm A coalgebra $C$ is called a left {\it Calabi-Yau coalgebra of dimension
$n$} (simply written as CY-$n$) if

(i) $\D^b_{fd}({}^C\!\M)$ is Hom-finite;

(ii) $\D^b_{fd}({}^C\!\M)$ is a Calabi-Yau category of dimension
$n$.}
\end{defn}

Similarly, we can define right Calabi-Yau coalgebras. Note that for a  semiperfect or an artinian coalgebra $C$,  $C$ is left Calabi-Yau if and only if it is right Calabi-Yau. But in general,  we don't know whether or not left Calab-Yau and right Calabi-Yau are equivalent. In the following, a CY coalgebra always means a left CY coalgebra.

We now list some basic properties of a Calabi-Yau coalgebra:

\begin{prop}
{\rm(i)} If $C$ and $D$ are Morita-Takeuchi equivalent, then $C$ is
CY-$n$ if and only if $D$ is CY-$n$;

{\rm(ii)} If $C$ is CY-$n$, then the global dimension $\text{\rm
gldim}C=n$;

{\rm(iii)} If $\{C_\lambda\}_{\lambda\in\Lambda}$ is a set of CY-$n$
coalgebras, then $\op_{\lambda\in\Lambda}C_\lambda$ is a CY-$n$
coalgebra.
\end{prop}
\proof the statement (i) is trivial. The statement (ii) follows from the fact that
the global dimension of $C$ is equal to the supremum of the injective
dimension of simples comodules (cf. \cite{NTZ}).

(iii) Let $C$ be a coalgebra, and $e_1,e_2\in C^*$ be a pair of
central orthogonal idempotents. We have that for any objects $X,Y\in
\D^+({}^C\!\M)$, $\Hom_{\D^+({}^C\!\M)}(Xe_1,Ye_2)=0$. Now let
$C=\op_{\lambda\in\Lambda}C_\lambda$, and let $e_\lambda\in C^*$ be the
central idempotent whose restriction to $C_\lambda$ is the counit
$\varepsilon_\lambda$ of $C_\lambda$ and $e_\lambda$ sends $C_\beta$ to zero
if $\beta\neq\lambda$. For $X\in \D^b_{fd}({}^C\!\M)$, note that
there are only finitely many idempotents
$e_{\lambda_1},\dots,e_{\lambda_n}$ such that $Xe_{\lambda_i}\neq0$
in $\D^b_{fd}({}^C\!\M)$. Hence we have the following natural isomorphisms for
any $X,Y\in\D^b_{fd}({}^C\!\M)$,
$$\begin{array}{ccl}
        \Hom_{\D^b_{fd}({}^C\!\M)}(X,Y)&=&\Hom_{\D^b_{fd}({}^C\!\M)}(\op_{\lambda\in\Lambda}Xe_\lambda,\op_{\lambda\in\Lambda}Ye_\lambda)\\
        &\cong&\op_{\lambda\in\Lambda}\Hom_{\D^b_{fd}({}^C\!\M)}(Xe_\lambda,Ye_\lambda)\\
        &\cong&\op_{\lambda\in\Lambda}\Hom_{\D^b_{fd}({}^C\!\M)}(Ye_\lambda,Xe_\lambda[n])^*\\
        &\cong&\Hom_{\D^b_{fd}({}^C\!\M)}(\op_{\lambda\in\Lambda}Ye_\lambda,\op_{\lambda\in\Lambda}Xe_\lambda[n])^*\\
        &\cong&\Hom_{\D^b_{fd}({}^C\!\M)}(Y,X[n])^*.
       \end{array}
$$
Therefore (iii) holds. \qed

\section{Calabi-Yau coalgebras of dimension 0 and 1}

In this section, we will see that the structures of the Calabi-Yau coalgebras of dimensions 0 and 1 are quite simple. A CY-$0$ coalgebra is nothing but a cosemisimple coalgebra as we will see in the next proposition.

\begin{prop} A coalgebra $C$ is CY-0 if and only if $C$ is
cosemisimple.
\end{prop}
\proof We only need to prove that a cosemisimple coalgebra is CY-0.
By \cite[Cor.3.6]{CDN}, $C$ is a symmetric coalgebra. Hence, for any
left $C$-comodule $N$, we have the left $C^*$-isomorphisms
$N^*\cong\Hom_{C^*}(N,C)\cong\Hom_{C^*}(N,C^*)$ (cf. \cite[Theorem
5.3]{CDN}). Since $C$ is cosemisimple, $\D^b_{fd}({}^C\!\M)$ is
exactly the homotopy category of bounded complexes with finite
dimensional cohomology, and any complex in $\D^b_{fd}({}^C\!\M)$ is
split. It is sufficient to show that the natural isomorphisms in
Definition \ref{cdef1} hold for all finite dimensional comodules.
Let $N$ and $M$ be any finite dimensional comodules. We have the following natural
isomorphisms $$\begin{array}{ccl}
                 \Hom_C(M,N)&\cong&\Hom_{C^*}(M,N)\\
                 &\cong&\Hom_{C^*}(M,N^{**})\\
                 &\cong&\Hom(M\ot_{C^*}N^*,\k)\\
                 &\cong&\Hom(M\ot_{C^*}\Hom_{C^*}(N,C^*),\k)\\
                 &\cong&\Hom_{C^*}(N,M)^*.
               \end{array}
$$ It follows that $C$ is CY-0.\qed

\begin{cor}\label{cor2.2} Let $C$ be a semiperfect coalgebra. Then $C$ is CY if
and only if $C$ is cosemisimple.
\end{cor}
\proof Since $C$ is semiperfect, the injective envelop of a simple
comodule is finite dimensional. Hence there are enough finite
dimensional injective comodules. Let $E$ be a finite dimensional
injective comodule. If $C$ is CY-$n$, then
$\Ext_C^n(E,E)\cong\Hom_C(E,E)^*\neq 0$. This implies that $n=0$. Therefore $C$ is cosemisimple. \qed

\begin{rem} From the proof of Corollary \ref{cor2.2}, we see that if
$C$ is a non-cosemisimple CY coalgebra, then any injective comodule
should be of infinite dimension.
\end{rem}

It is well-known that a hereditary coalgebra over an algebraically closed field is Morita-Takeuchi equivalent to a path coalgebra (cf. \cite{Chin}), a CY-1 coalgebra over an algebraically closed field must be Morita-Takeuchi equivalent to a path coalgebra. Let us recall some notations about quivers and path coalgebras. Let $Q$ be a quiver with the set of vertices $Q_0$ and the set of arrows $Q_1$. For an arrow $a\in Q_1$, we use $s(a)$ to denote the source of $a$ and $t(a)$ to denote the target of $a$. A nontrivial path is a sequence of arrows $p=a_1a_2\cdots a_n$ with $s(a_{i+1})=t(a_i)$. We say that the length of $p$ is $n$. Sometimes a vertex is called a path of length 0. For a path $p=a_1a_2\cdots a_n$, we define $s(p)=s(a_1)$ and
$t(p)=t(a_n)$. Denote by $\k Q$ the path algebra of $Q$ (in general, $\k Q$ has no unit). The the multiplication of two pathes $p,q$ is defined as $pq$ if $s(q)=t(p)$ and $0$ otherwise.
Let $a$ be an arrow and $p=a_1a_2\cdots a_n$
be a path. We define $pa^{-1}=a_1\cdots a_{n-1}$ if $a=a_n$ and $0$ otherwise. If $x=k_1p_1+\cdots+k_np_n$ is a linear combination of paths, then we define $xa^{-1}=k_1p_1a^{-1}+\cdots+k_np_na^{-1}$.
Similarly, if $q$ is a general path, we define $pq^{-1}=r$ if $p=rq$.

Following the notations of \cite{JMN}, we use $CQ$ to denote the
path coalgebra of $Q$. There is a nondegenerated bilinear map
$\langle\ ,\ \rangle:CQ\times\k Q\longrightarrow\k$ defined by
$\langle p,q\rangle=\delta_{p,q}$ (the Kronecker delta), where $p,q$
are two pathes of the quiver $Q$. The bilinear form $\langle\ ,\
\rangle$ induces an injective map $\iota:\k Q\longrightarrow
(CQ)^*$. Clearly, $\iota$ preserves the multiplications (cf.
\cite{JMN}). Since $CQ$ is a $(CQ)^*$-bimodule, for $x\in CQ$ and
$y\in \k Q$, the notion $\iota(y)x$ (or $x\iota(y)$) makes sense. We
list two simple properties of the bilinear form $\langle\ ,\
\rangle$, which will be frequently used.

\begin{lem} {\rm(i)} Let $a$ be an arrow, $x\in CQ$ and
$y\in \k Q$. Then $\langle xa^{-1},y\rangle=\langle x,ya\rangle$.

{\rm(ii)} Let $x\in CQ$ and $y\in \k Q$. Assume
$y=k_1p_1+\cdots+k_np_n$ with $k_1k_2\cdots k_n\neq0$ and
$s(p_1)=\cdots=s(p_n)$. If $\iota(y)x=0$ then $\langle
x,y\rangle=0$.
\end{lem}
\proof (i) is obvious. The statement (ii) follows from (i)
and the fact that $xp^{-1}=\iota(p)x$ for any path $p$. \qed

Let $i\in Q_0$ be a vertex, let $e_i\in (CQ)^*$ be the idempotent
corresponding to $i$, and let $S_i$ be the simple left $CQ$-comodule
corresponding to $i$. For simplicity, we write $C$ for $CQ$. The minimal injective resolution of $S_i$ can be written as follows:
\begin{equation}\label{cequ1}
    0\longrightarrow S_i\longrightarrow e_iC\overset{f}\longrightarrow
\bigoplus_{a\in Q_1, t(a)=i}e_{s(a)}C\to 0,
\end{equation}
where $f$ is defined by $\displaystyle f(x)=\sum_{t(a)=i}\iota(a)x$
in which the multiplication is the left $C^*$-module action, and we
regard $\iota(a)x$ as an element in $e_{s(a)}C$. Here we need to point out that if two arrows $a$ and $b$ have the same source and the same target in the above sequence, we should distinguish $e_{s(a)}C$ from $e_{s(b)}C$ in the direct sum $\bigoplus_{a\in Q_1;
t(a)=i}e_{s(a)}C$.

\begin{prop} Let $C$ be a coalgebra. Then $C$ is CY-1 if and only if $C$ is Morita-Takeuchi equivalent to a direct sum of copies of $\k[x]$, where $\k[x]$ is the path coalgebra of the quiver $Q$ with one vertex and one arrow.
\end{prop}
\proof Note that $\k[x]$ is a CY-1 coalgebra because the dual
algebra of $\k[x]$ is the power series algebra which is noetherian
and CY-1 (cf. \cite{HTVZ}). Hence a direct sum of copies
of $\k[x]$ is also CY-1.

Conversely, assume that $C$ is CY-1. Without loss of generality, we
can further assume that $C$ is the path coalgebra of a quiver $Q$. Let $i,j\in Q_0$. By the CY property,
$\dim\Ext_C^1(S_j,S_i)=\dim\Ext_C^0(S_i,S_j)=\delta_{i,j}$. From the
minimal resolution (\ref{cequ1}) of $S_i$, we have
$\dim\Ext^1_C(S_j,S_i)=\#\{a\in Q_1|s(a)=j,t(a)=i\}$. It follows that for vertices $i,j$ and $i\neq j$, there exist no
arrows from $i$ to $j$, and for each vertex $i$ there is a unique
arrow from $i$ to itself. Therefore, $C$ is a direct sum of copies of
$\k[x]$. \qed

\section{Minimal injective resolutions of the simple comodules}

In order to investigate CY-2 and CY-3 coalgebras, we need more notations about path
coalgebras of quivers with relations. Recall from \cite{Sim} that a
{\it quiver with relations} means a pair $(Q,\Omega)$ where $\Omega$
is a two-side ideal of $\k Q$ contained in $\k Q_{\ge2}$. Note that
the path algebra $\k Q$ is graded. If $\Omega$ is a graded ideal,
then we say that the relations are {\it homogeneous}. The path coalgebra
of $(Q,\Omega)$ is the following subcoalgebra of $CQ$:
$$C(Q,\Omega)=\{x\in CQ|\langle x,\Omega\rangle=0\}.$$
It is well known that any coalgebra is Morita-Takeuchi equivalent to
a basic coalgebra and a basic coalgebra over an algebraically closed
field is isomorphic to a subcoalgebra of the path coalgebra of its
Gabriel quiver (cf. \cite{CM}). However not every subcoalgebra $C$
of a path coalgebra is of the form $C=C(Q,\Omega)$ (cf. \cite{JMN}).
We wish to understand when a path coalgebra of a quiver with
relations is CY-2 or CY-3. To this end, we need to investigate the
minimal injective resolutions of simple comodules. For convenience, we introduce some temporary notations. Let
$p=a_1a_2\cdots a_n$ be a nontrivial path. We define
$lead(p)=\{a_n\}$, the set of leading arrow. If $p$ is a trivial path, then we let $lead(p)=\emptyset$ . For $x\in\k Q$ and
$x=k_1p_1+\cdots+k_np_n$ with $k_1k_2\cdots k_n\neq0$, we define $lead(x)=lead(p_1)\cup\dots\cup
lead(p_n)$. If $S$ is a subset of $\k Q$, then we define
$lead(S)=\bigcup_{x\in S}lead(x)$.

We choose a generating set of the relation ideal $\Omega$ in the
following way. Denote by $\k Q_{\leq n}$ the subspace of
linear combinations of paths of length less or equal to $n$. Then
the path algebra $\k Q$ is a filtered algebra. The restriction of
the filtration to $\Omega$ results in a filtration on $\Omega$. Set
$\Omega(n)=\Omega\cap\k Q_{\leq n}$. Note that
$\Omega(0)=\Omega(1)=0$. For $i,j\in Q_0$, let $\k Q_{i,j}$ be the
subspace of all the linear combinations of paths from $i$ to $j$.
Set $\Omega(n)_{i,j}=\Omega(n)\cap\k Q_{i,j}$. Then
$\Omega(n)=\bigoplus_{i,j\in Q_0}\Omega(n)_{i,j}$. Choose a basis
$B(2)_{i,j}$ of $\Omega(2)\cap\k Q_{i,j}$ for all $i,j\in Q_0$. Set
$R(2)=\bigcup_{i,j\in Q_0}B(2)_{i,j}$. Then $R(2)$ is a basis of
$\Omega(2)$. For $n\ge2$, let $I(n)$ be the ideal of $\k Q$
generated by $\Omega(n)$. Set $I(n)_{i,j}=I(n)\cap\k Q_{i,j}$. Then
$I(n)=\bigoplus_{i,j\in Q_0}I(n)_{i,j}$. For $i,j\in Q_0$, let
$V(n+1)_{i,j}$ be a subspace of $\Omega(n+1)_{i,j}$ such that
$\Omega(n+1)_{i,j}=(I(n)_{i,j}\cap\Omega(n+1)_{i,j})\op
V(n+1)_{i,j}$. Choose a basis $B(n+1)_{i,j}$ of $V(n+1)_{i,j}$. Set
$R(n+1)=\bigcup_{i,j\in Q_0}B(n+1)_{i,j}$. Now let
$R=\bigcup_{n\ge2}R(n)$. The set $R$ possesses the following properties:
\begin{enumerate}
\item[(i)] $R$ generates the ideal $\Omega$, and $R$ is minimal, that is,
any proper subset of $R$ can not generates $\Omega$;
\item[(ii)] each element of $R$ is a combination of paths with common
source and common target.
\end{enumerate}

We call a subset $R$ of $\Omega$ satisfying the properties (i) and (ii)  {\it a minimal set of relations} of $(Q,\Omega)$.

Because of the property (ii) above, for an element $r\in R$ we may write
$s(r)$ for the common source and $t(r)$ for the common target. Note
that if $\Omega$ is a graded ideal then we may choose the minimal set $R$ to be homogeneous, that is, each element of $R$ is a combination of
paths with the same length.

Let $(Q,\Omega)$ be a quiver with relations. We say that the relation ideal $\Omega$ is {\it locally finite} if there is a minimal set of relations $R$ such that for every pair of vertices $(i,j)$ the set $\{r\in
R|s(r)=i,t(r)=j\}$ is a finite set or empty. Note that the locally
finite property of $\Omega$ is independent of the choice of the
minimal set of relations $R$.

Now we are ready to construct the minimal injective resolution of
a simple comodule of certain path coalgebra. Let
$(Q,\Omega)$ be a quiver with relations, and $C=C(Q,\Omega)$. Let $R$ be a minimal set of
relations. Assume that $\Omega$ is locally finite. Let $S_n$ be the simple comodule
corresponding to the vertex $n$. Let $R_n=\{r\in R|t(r)=n\}$. We
construct a sequence:
\begin{equation}\label{cequ2}
    0\longrightarrow S_n\longrightarrow e_nC\overset{f}\longrightarrow \bigoplus_{a\in
    Q_1,t(a)=n}e_{s(a)}C\quad\overset{g}\longrightarrow \bigoplus_{r\in
    R_n}e_{s(r)}C,
\end{equation}
where the first map is the embedding map, and $f$ and $g$ are
constructed as follows. For $x\in e_nC$, $\displaystyle
f(x)=\sum_{t(a)=n,a\in Q_1}\iota(a)x$, where $\iota(a)x$ is regarded as an
element in $e_{s(a)} C$. Similar to the sequence (\ref{cequ1}), we
should distinguish $e_{s(a)} C$ from $e_{s(b)} C$ if $a$ and $b$
have the same source and the same target. For $x\in e_{s(a)}C$ with
$a\in\hspace{-9pt}/\, lead(R_n)$, define $g(x)=0$. If $x\in e_{s(a)}C$ with $a\in lead(R_n)$, then we define $g(x)=\sum_{r\in R_n}\iota(ra^{-1})x$.
Since $\Omega$ is locally finite, $g$ is well defined. Here
$\iota(ra^{-1})x$ is viewed as an element in $e_{s(r)}C$, and we also need to distinguish $e_{s(r)}C$ from $e_{s(r')}C$ in the direct sum
$\bigoplus_{t(r)=n}e_{s(r)}C$ if $r$ and $r'$ share the same source (and the same target).

\begin{thm}\label{thm3.1} Let $(Q,\Omega)$ be a quiver with relations such that
$\Omega$ is locally finite. Then the sequence (\ref{cequ2})
constructed above is exact.
\end{thm}
\proof We first prove $gf=0$. For $x\in e_nC$, we may write $x$ as
$x=\sum_{i=1}^sk_ip_i+\sum_{j=1}^mt_jq_j$, where $p_i$'s are paths
such that $lead(p_i)\nsubseteq lead(R_n)$ for all $i=1,\dots,s$ and
$q_j$'s are paths such that $lead(q_j)\subseteq lead(R_n)$ for all
$j=1,\dots,m$. Let $a_j$ be the leading arrow of $q_j$. We have
$$\begin{array}{ccl}
    gf(x)&=&gf(\sum_{j=1}^mt_jq_j)\\
    &=&g(\sum_{j=1}^mt_j\iota(a_j)q_j)\\
    &=&\sum_{j=1}^mt_j(\sum_{r\in
    R_n}\iota(ra_j^{-1})\iota(a_j)q_j)\\
    &\overset{(I)}=&\sum_{j=1}^mt_j(\sum_{r\in
    R_n}\iota(r)q_j)\\
    &=&\sum_{r\in
    R_n}\iota(r)x\\
    &=&\sum_{r\in
    R_n}\sum_{(x)}x_{(1)}\langle x_{(2)},r\rangle\\
    &=&0,
  \end{array}
$$ where the last identity holds because $C$ is a coalgebra of $\k
Q$ and any element in $C$ is orthogonal with $\Omega$. For the identity $(I)$ one may check it straightforward. We should point out that $\iota(ra^{-1}_j)\iota(a_j)$ may not equal to $\iota(r)$.

Next we show the inclusion $\ker g\subseteq\text{im} f$. Given an element $y\in \ker g$, then $y$ is in one of the following two cases.

Case (1): if $y\in e_{s(a)}C$ with $a\in\hspace{-8pt}/ \, lead(R_n)$, we let $x=ya$. In this case it is clear that $\iota(r)x=0$ for all $r\in R$.

Case (2): if $y\in \bigoplus_{a\in lead(R_n)}e_{s(a)}C$, we can assume that $y=k_1\alpha_1+\cdots+k_m\alpha_m$ with
$k_1k_2\cdots k_m\neq0$ and $\alpha_i\in e_{s(a_i)}C$
($i=1,\dots,m$). In this case we have
$$\begin{array}{ccl}
   g(y)&=&g(\sum_{i=1}^mk_i\alpha_i)\\
   &=&\sum_{i=1}^mk_i(\sum_{r\in
   R_n}\iota(ra_i^{-1})\alpha_i)\\
   &=&\sum_{i=1}^mk_i\left(\sum_{r\in
   R_n}\iota(ra_i^{-1})\iota(a_i)(\alpha_ia_i)\right)\\
   &=&\sum_{i=1}^mk_i\sum_{r\in
   R_n}\iota(r)(\alpha_ia_i)\\
   &=&\sum_{r\in
   R_n}\iota(r)\sum_{i=1}^m(k_i\alpha_ia_i),
  \end{array}
$$
where the forth identity holds because $\alpha_ia_i$ is a linear
combination of paths with the same leading arrow $a_i$. Let
$x=\sum_{i=1}^mk_i\alpha_ia_i$. Recall that $\iota(r)x\in e_{s(r)}C$. Hence $g(y)=0$ implies $\iota(r)x=0$ for all
$r\in R_n$. Of course, if $r\in\hspace{-9pt}/ \,  R_n$, we certainly have
$\iota(r)x=0$.

To show that in both cases $y\in \mathrm{im}f$, we only need to show that $x\in e_nC$ since $y=f(x)$.  This is equivalent to proving that $\langle x,\Omega\rangle=0$. Indeed, for any $\beta\in \Omega$, we write
$\beta=\sum_{i=1}^sk_ip_ir_i+\sum_{j=1}^mt_jp_j'r_j'q_j$ where
$r_i,r'_j\in R$ and $q_j$ is a nontrivial path for every $j$. Let
$a_j$ be the leading arrow of $q_j$ for $j=1,\dots,m$. We have
$$\begin{array}{ccl}
\langle x,\beta\rangle&=&\displaystyle\sum_{i=1}^sk_i\langle
x,p_ir_i\rangle+\sum_{j=1}^mt_j\langle
x,p_j'r_j'q_j\rangle\\
&=&\displaystyle\sum_{i=1}^sk_i\langle
\iota(r_i)x,p_i\rangle+\sum_{j=1}^mt_j\langle
\iota(a_j)x,p_j'r_j'q_ja_j^{-1}\rangle\\
&=&0,
  \end{array}
$$ where the last identity holds because $\iota(a_j)x\in C$ and
$p_j'r_j'q_ja_j^{-1}\in \Omega$. This completes the proof. \qed

\begin{rem} There is a right version of the sequence (\ref{cequ2}).
Let $C=C(Q,\Omega)$ be as in Theorem \ref{thm3.1}, and let $S_m$ be
the right simple comodule corresponding to the vertex $m$. The first
two steps of the minimal injective resolution of $S_m$ is:
\begin{equation}\label{cequ3}
 0\longrightarrow S_m\longrightarrow Ce_m\overset{f}\longrightarrow \bigoplus_{a\in
    Q_1,s(a)=n}Ce_{t(a)}\quad\overset{g}\longrightarrow \bigoplus_{r\in
    R,s(r)=m}Ce_{t(r)},
\end{equation}
where the maps $f$ and $g$ are defined as follows.
Define $\displaystyle f(x)=\sum_{s(a)=m, a\in Q_1}x\iota(a)$, for $x\in Ce_m$, and for $x\in
Ce_{t(a)}$, define $\displaystyle
g(x)=\sum_{s(r)=m,r\in R}x\iota(a^{-1}r)$.
\end{rem}

Since all the items (except the simple comodule $S_n$) in the
sequence (\ref{cequ2}) are injective and the socle of each injective
comodule is contained in the image of the map, the sequence
(\ref{cequ2}) is the first two steps of the minimal injective
resolution of $S_n$.

\begin{cor} Let $(Q,\Omega)$ be a quiver with relations. If $C=C(Q,\Omega)$ is CY, then, for any pair of vertices
$(i,j)$, there are at most finitely many arrows from $i$ to $j$.
\end{cor}
\proof The sequence $$0\longrightarrow S_n\longrightarrow
e_nC\overset{f}\longrightarrow \bigoplus_{a\in
Q_0,t(a)=n}e_{s(a)}C$$ as a part of the sequence (\ref{cequ2}) is
always exact for any path coalgebra of a quiver with relations. Then
the result follows from the hypothesis that the derived category of
complexes with finite dimensional cohomology is Hom-finite.\qed

\section{Calabi-Yau coalgebras of dimensions 2 and 3.}

With the preparation of the preceding sections, we can now deduce some necessity conditions for a quiver with
relations so that its path coalgebra is CY-2 or CY-3.  The following theorems are dual to the corresponding results
in \cite{Boc}.

\begin{thm} \label{cthm1} Let $(Q,\Omega)$ be a quiver with locally finite
relation ideal $\Omega$. Let $R$ be a minimal set of relations of
$(Q,\Omega)$. Assume that $C=C(Q,\Omega)$ is CY-2. Then we have the following.
\begin{enumerate}
\item[(i)] For each vertex $n$, there is a unique element $r\in R$
such that $s(r)=t(r)=n$;
\item[(ii)] let $W=\{ra^{-1}|a\in Q_1,r\in R\}$, and let
$\overline{\Omega}=(W)$ be the ideal generated by the elements in
$W$. Then $\overline{\Omega}=\k Q_{\ge1}$;
\item[(iii)] For any two vertices $n,m$, $\#\{a\in Q|s(a)=n, t(a)=m\}=\#\{b\in Q_1|s(b)=m, t(b)=n\}$; for any vertex $n$, there are at most finitely many arrows starting from $n$, and at most finitely many
arrows ending at $n$.
\item[(iv)] If $\Omega$ is a graded ideal of $\k Q$, then any element
$r$ in $R$ is a linear combination of paths of length 2.
\end{enumerate}
\end{thm}
\proof Since $C$ is of global dimension 2, the minimal injective
resolution of $S_n$ is:
$$0\longrightarrow S_n\longrightarrow e_nC\overset{f}\longrightarrow \bigoplus_{a\in
    Q_1,t(a)=n}e_{s(a)}C\quad\overset{g}\longrightarrow \bigoplus_{r\in
    R_n}e_{s(r)}C\longrightarrow0.$$ Now by the CY property, we have
$\dim\Ext_C^2(S_m,S_n)=\dim\Hom_C(S_n,S_m)=\delta_{n,m}$. Hence (i)
follows.

(ii) It is not hard to see $\overline{\Omega}\supseteq \Omega$.
Hence the coalgebra  $\overline{C}$ defined by $\{x\in CQ|\langle x,\overline{\Omega}\rangle=0\}$ is a subcoalgebra of
$C$. View $\overline{C}$ as a left $C$-comodule. Following (i), the
minimal injective resolution of $S_n$ reads as follows:
$$0\longrightarrow S_n\longrightarrow e_nC\overset{f}\longrightarrow \bigoplus_{a\in
    Q_1,t(a)=n}e_{s(a)}C\quad\overset{g}\longrightarrow e_nC\longrightarrow0.$$
Let $M$ be a finite dimensional $C$-subcomodule of $\overline{C}$.
Applying $\Hom_C(M,-)$ to the above injective resolution of
$S_n$, we obtain that $\Ext^2_C(M,S_n)$ is the cokernel of
$g_*=\Hom_C(M,g)$. Assume that $r$ is the unique element in $R$ such that
$s(r)=t(r)=n$. We claim that $g_*=0$. Indeed, let $h\in\Hom_C(M,\bigoplus_{a\in
    Q_1,t(a)=n}e_{s(a)}C)$. Since $M$ is finite dimensional,
the image of $h$ lies in $\op_{i=1}^ke_{s(a_i)}C$ for finitely many
arrows $a_1,\dots,a_n$. We still use $h$ to denote the induced
morphism $h:M\longrightarrow \op_{i=1}^ke_{s(a_i)}C$. By Prop.
\ref{aprop} in the appendix, we have the following diagram
$$\xymatrix{
  \Hom_C(M,\op_{i=1}^ke_{s(a_i)}C) \ar[d]_{\cong} \ar[r]^{\qquad g_*} & \Hom_C(M,e_nC) \ar[d]^{\cong} \\
  h_C(\op_{i=1}^ke_{s(a_i)}C,M)^*\ar[d]_{\cong}\ar[r]^{\qquad
  h_C(g,M)^*}&h_C(e_nC,M)^*\ar[d]^{\cong}\\
  \op_{i=1}^k(Me_{s(a_i)})^* \ar[r]^{\quad\theta} & (Me_n)^*,   }$$
where $\theta$ is the dual of the morphism: $\xymatrix@C=1.5cm{ Me_n
\ar[r]^{\cdot\iota(ra_i^{-1})\qquad} &  \op_{i=1}^kMe_{s(a_i)}. }$
Recall that $M$ is contained in $\overline{C}$, and so is $Me_n$.
Hence $\langle Me_n,\overline{\Omega}\rangle=0$, and
$Me_n\iota(ra_i^{-1})=0$ for all $i$. Thus the map $\theta$ in the
diagram is the zero map. It follows that the map $g_*$ is  the zero
map as well. So the claim follows. Thus we have
$\Ext^2_C(M,S_n)\cong\Hom_C(M,e_nC)$. By the CY property,
$\Ext^2_C(M,S_n)\cong\Hom_C(S_n,M)^*$. Therefore
\begin{equation}\label{cequ4}
1\ge\dim\Hom_C(S_n,M)=\dim\Hom(M,e_nC)=\dim h_C(e_nC,M)=\dim Me_n.
\end{equation}
Suppose that there is a nonzero element
$\overline{c}\in\overline{C}$ such that it is a combination of nontrivial
paths, say, $\overline{c}=k_1p_1+\cdots+k_mp_m$ with $k_1\neq0$.
Let $D$ be the subcoalgebra of $\overline{C}$ generated by
$\overline{c}$. Then $D$ can be viewed as a finite dimensional left
subcomodule of $\overline{C}$. Since $D$ must contain the vertex
$s(p_1)$ and $\overline{c}$, we deduce that $\dim De_{s(p_1)}$ is at
least 2. This contradicts with the fact (\ref{cequ4}). Hence
$\overline{C}$ is exactly the subcoalgebra of $CQ$ generated by
the vertices of $Q$. Therefore $\overline{\Omega}$ must be the ideal
of $\k Q$ generated by all the arrows of $Q$.

(iii) The first part follows from the fact $\dim \Ext^1_C(S_n,S_m)=\dim\Ext^1_C(S_m,S_n)$. The second part follows from (i)
and (ii) since the arrows in $\overline{\Omega}$ starting from
a vertex $n$ are contained in the ideal generated by
$\{ra^{-1}|a\in Q_1\}$, where $r$ is the unique element $r$ in $R$
with $s(r)=n$.

(iv) If $\Omega$ is graded, then any element in $R$ is homogeneous.
For any element $r\in R$, we may assume $r=k_1p_1+\cdots+k_mp_m$ with
$k_1k_2\cdots k_m\neq0$. By (i) and (ii), there are at least one
path among $p_1,\dots,p_m$ of length 2. Otherwise
$\overline{\Omega}$ could not be $\k Q_{\ge1}$. This forces all the
paths $p_1,\dots,p_m$ to be of length 2. \qed

\begin{thm} Let $(Q,\Omega)$ be as in Theorem \ref{cthm1}.
Assume $C=C(Q,\Omega)$ is CY-3.

{\rm(i)} Let $R$ be a minimal relation set. For any vertices $i,j$,
we have $\#\{a\in Q_1|s(a)=i,t(a)=j\}=\#\{r\in R|s(r)=j,t(r)=i\}$;

{\rm(ii)} Assume further that the ideal $\Omega$ is graded, and for any vertices $i,j\in Q_0$ and any integer $n\ge1$, there are only finitely many paths of length $n$ starting from $i$ and ending at $j$.
Then we may choose a minimal relation set $R$ such that every
element of $R$ is a combination of paths with a fixed length;

{\rm(iii)} Under the assumptions of (ii), we may choose a minimal
relation set $R$ and a correspondence $\nu:Q_1\to R$ such that
$s(\nu(a))=t(a), t(\nu(a))=s(a)$ for all $a\in Q_1$, and for each
arrow $b$ with $t(b)=i$, $\displaystyle
r_b=\sum_{s(a)=i,a\in Q_1}k_b ar_ab^{-1}$, where $r_a=\nu(a)$ and $k_b\in \k$.

\end{thm}
\proof (i) By the CY property,
$\dim\Ext_C^1(S_j,S_i)=\dim\Ext^2_C(S_i,S_j)$. From the minimal
injective resolutions of $S_i$ and of $S_j$, one easily obtain
$\#\{a\in Q_1|s(a)=i,t(a)=j\}=\dim\Ext_C^1(S_j,S_i)$ and $\#\{r\in
R|s(r)=j,t(r)=i\}=\dim\Ext^2_C(S_i,S_j)$.

(ii) Since $\Omega$ is graded, we may choose a minimal relation set
$R$ such that every element of $R$ is a combination of paths with
the same length. Since the global dimension of $C$ is 3, the minimal
injective resolution of $S_i$ is of the following form
\begin{equation}\label{cequ5}
0\longrightarrow S_i\longrightarrow e_iC\overset{f}\longrightarrow
\bigoplus_{t(a)=i}e_{s(a)}C\quad\overset{g}\longrightarrow
\bigoplus_{r\in
    R,t(r)=i}e_{s(r)}C\overset{\eta}\longrightarrow e_iC\longrightarrow0,
\end{equation}
    where $f$ and $g$ is the map formed in (\ref{cequ2}). We want to
construct the map $\eta$ explicitly. Since $C$ is graded, all the maps
in the above sequence are graded maps. The map $p$ is determined by
a sequence of maps $\eta_r:e_{s(r)}C\to e_iC$ for $r\in R$ and
$t(r)=i$. Note that $\Hom_C(e_{s(r)}C,e_iC)\cong
h_C(e_iC,e_{s(r)}C)^*\cong (e_{s(r)}Ce_i)^*$. Since $\eta$ is a graded
map, the finiteness assumption on $Q$ implies that there is homogeneous element
$\alpha_{r}\in e_{s(r)}Ce_i$ such that $\eta_r(x)=\iota(\alpha_{r})x$
for all $x\in e_{s(r)}C$, where we view $e_{s(r)}Ce_i$ as a subset
of $\k Q$. Let $\overline{\Omega}$ be the ideal of $\k Q$ generated
by the set $R\bigcup\{\alpha_{r}|r\in R\}$, and let
$\overline{C}=\{x|\langle x,\overline{\Omega}\rangle=0\}$. Then $\overline{C}$ is a subcoalgebra of $C$. For any finite dimensional subcoalgebra $D$ of
$\overline{C}$, a similar argument to the one in the proof of Theorem
\ref{cthm1} shows that $\Ext_C^3(D,S_i)\cong\Hom_C(D,e_iC)$. By the CY
property we obtain $$\dim De_i=\dim
h_c(e_iC,D)=\dim\Hom_C(D,e_iC)=\dim\Hom_C(S_i,D)\leq 1$$ for all
$i\in Q_0$, which implies that $\overline{C}$ must be the coradical
of $C$. Hence $\overline{\Omega}=\k Q_{\ge1}$. Moreover, for any
vertices $i,j$, we have $\#\{\alpha_r|s(r)=j,t(r)=i\}=\#\{a\in
Q_1|s(a)=i,t(a)=j\}$ by (i). Since $\Omega\subseteq \k Q_{\ge2}$, we obtain that
$\alpha_r$ is a linear combination of arrows starting from $t(r)$ and
ending at $s(r)$ for all $\alpha_r$, and that the elements in
$\{\alpha_r|s(r)=j,t(r)=i\}$ are linearly independent. Since all the maps in the resolution (\ref{cequ5}) are graded the
statement (ii) follows.

(iii) Following (ii), we may choose a minimal relation set $R$
whose elements are linear combination of paths with fixed length.
Also the map $p_r:e_{s(r)}C\to e_iC$ is determined by an element that is a
linear combination of arrows. Once again, for any pair of vertices
$(i,j)$, the facts that $\#\{\alpha_r|s(r)=j,t(r)=i\}=\#\{a\in
Q_1|s(a)=i,t(a)=j\}$ and that the elements in
$\{\alpha_r|s(r)=j,t(r)=i\}$ are linearly independent enable us to
combine the elements in $R_{j,i}:=\{r\in R|s(r)=j,t(r)=i\}$ linearly
to obtain a new set $R'_{j,i}$ such that $\eta_{r'}:e_{s(r')}C\to e_iC$
is defined by $\eta_{r'}(x)=\iota(a)x$ with $a$ an arrow starting from
$j$ and ending at $i$ for all $r'\in R'_{j,i}$. Moreover, we have established
a correspondence $\nu_{i,j}:\{a\in Q_1|s(a)=i,t(a)=j\}\to R'_{j,i}$.
Now let $R'=\bigcup_{i,j} R'_{i,j}$. We obtain a correspondence
$\nu:Q_1\to R'$ such that $\eta_{\nu(a)}:e_{s(\nu(a))}C\to e_iC$ is
defined by $\eta_{\nu(a)}(x)=\iota(a)x$. For simplicity, we assume that
$R$ itself has the above properties, and write $r_a=\nu(a)$ for all
$a\in Q_1$. Let $a\in Q_1$ with $t(a)=i$. For any $x\in e_{s(a)}C$,
we have
$$0=\eta g(x)=\eta\left(\sum_{s(d)=i}\iota(r_da^{-1})x\right)
=\sum_{s(d)=i}\iota(d)\iota(r_da^{-1})x=\sum_{s(d)=i}\iota(dr_da^{-1})x.$$
This implies $\langle e_{s(a)}C,\sum_{s(d)=i}dr_da^{-1}\rangle=0$. Since
$\sum_{s(d)=i}dr_da^{-1}$ is a linear combination of paths with the
same target $s(a)$, we obtain $\langle
C,\sum_{s(d)=i}dr_da^{-1}\rangle=0$. Hence
$\sum_{s(d)=i}dr_da^{-1}\in \Omega$ by the finiteness assumption on $Q$.

For $\xi\in \Ext^2_C(S_i,S_j)$ and $\gamma\in \Ext_C^1(S_j,S_i)$, we
compute the Yoneda products $\xi*\gamma$ and $\gamma*\xi$. Assume
there are $n$ arrows from $j$ to $i$ labeled as $a_1,\dots,a_n$, and
$m$ arrows starting from $i$ labeled as $d_1,\dots,d_m$. Note that
$\dim\Ext^2_C(S_i,S_j)=n=\dim\Ext_C^1(S_j,S_i)$. As we have seen that for every $u$ $(1\leq u\leq n)$,
$\sum_{v=1}^md_vr_{d_v}a_u^{-1}\in \Omega$, we obtain
\begin{equation}\label{cequ6}
    \sum_{v=1}^md_vr_{d_v}a_u^{-1}=\sum_{w=1}^nl_{uw}r_{a_w},\ \text{where
    $l_{uw}\in\k$}.
\end{equation}
We rewrite the minimal injective resolution of $S_i$ as follows:
{\small$$0\longrightarrow S_i\longrightarrow
e_iC\overset{f_i}\longrightarrow
\bigoplus_{u=1}^ne_{s(a_u)}C\quad\op\bigoplus_{t(a)=i,,a\neq
a_u}e_{s(a)}C\quad\overset{g_i}\longrightarrow
\bigoplus_{v=1}^me_{s(r_{d_v})}C\overset{\eta_i}\longrightarrow
e_iC\longrightarrow0.$$} Since
$\Ext_C^1(S_j,S_i)\cong\Hom_C(S_j,\bigoplus_{u=1}^ne_{s(a_u)}C)$,
$\gamma$ can be represented by a map
$f_\gamma\in\Hom_C(S_j,\bigoplus_{u=1}^ne_{s(a_u)}C)$.
similarly, $\xi$ can be represented by a map
$g_\xi\in\Hom_C(S_i,\bigoplus_{v=1}^me_{s(r_{d_v})}C)$. Consider now the
following diagram {\small$$\xymatrix{
   0 \ar[r] & S_j \ar[dr]^{f_\gamma}\ar[r]&e_jC\ar[d]^\theta\ar[r]^{f_j}&\displaystyle\bigoplus_{t(b)=j}e_{s(b)}C\ar[d]^\zeta\ar[r]^{g_j}&
   \displaystyle\bigoplus_{t(r)=j}e_{s(r)}C\ar[d]^\varphi\\
  S_i\ar[r]&e_iC \ar[r]_{f_i\qquad\qquad\quad\qquad} &\displaystyle\bigoplus_{u=1}^ne_{s(a_u)}C\quad\op\bigoplus_{t(a)=i,,a\neq
a_u}e_{s(a)}C\ar[r]_{\qquad\qquad g_i}&
\displaystyle\bigoplus_{v=1}^me_{s(r_{d_v})}C\ar[r]_{\quad
\eta_i}&e_iC. }$$} Now assume $f_\gamma:S_j\to
\bigoplus_{u=1}^ne_{s(a_u)}C$ is defined by $x\mapsto
(k_1x,\dots,k_nx)$ where  $k_ux$ is regarded as an element in
$e_{s(e_u)}C$. Now we can construct the maps $\theta$, $\zeta$ and
$\varphi$ as follows. For $x\in e_jC$, let
$\theta(x)=(k_1x,\dots,k_nx)$; for $x\in e_{s(b)}C$, let
$\zeta(x)=\sum_{u,v}k_u\iota(r_{d_v}a_u^{-1}b^{-1})x$. For the map $\varphi$,
we notice that
$\displaystyle\bigoplus_{t(r)=j}e_{s(r)}C=\bigoplus_{t(r)=j,s(r)\neq
i}e_{s(r)}C\op\bigoplus_{u=1}^ne_{s(r_{a_u})}C$. Now if $s(r)\neq
i$ and $x\in e_{s(r)}C$, set $\varphi(x)=0$, and if $x\in
e_{s(r_{a_w})}C$, set $\varphi(x)=\sum_{u=1}^nk_ul_{uw}x$, where
$l_{uw}$'s are the coefficients in the identity (\ref{cequ6}). Now
it is straightforward to check that the diagram above is
commutative. Suppose that $g_\xi\in
\Hom_C(S_i,\bigoplus_{v=1}^me_{s(r_{d_v})}C)$ is defined by
$g_\xi(x)=(\overline{k}_1x,\dots,\overline{k}_mx)$ where we view
$\overline{k}_vx$ as an element of $e_{s(r_{d_v})}C$. Then
$\gamma*\xi\in \Ext^3_C(S_i,S_i)$ is represented by the map
$\Psi:=\varphi\circ g_\xi\in\Hom_C(S_i,e_iC)$. Further, if $x\in S_i$, we
see that $\Psi(x)=\sum_{u,w=1}^n\overline{k}_wk_ul_{uw}x$.

Similarly, we can see that $\xi*\gamma$ is represented by a map
$\psi\in \Hom_C(S_j,e_jC)$ with
$\psi(y)=\sum_{u=1}^n\overline{k}_uk_uy$ for all $y\in e_jC$. Since
$C$ is CY-3, by \cite[Appendix]{Boc} there are trace map
$\text{Tr}_i:\Ext_C^3(S_i,S_i)\to \k$ for all $i$ such that
$\text{Tr}_i(\gamma*\xi)=\text{Tr}_j(\xi*\gamma)$. Since
$\Ext^3_C(S_i,S_i)$ is of dimension 1 for all $i$, the trace maps
$\text{Tr}_i:\Ext_C^3(S_i,S_i)\to \k$ is represented by a scalar
$\lambda_i$. Hence we have
$\lambda_i\sum_{u,w=1}^n\overline{k}_wk_ul_{uw}=\lambda_j\sum_{u=1}^n\overline{k}_uk_u$
for arbitrary choices of $(k_1,\dots,k_n)$ and
$(\overline{k}_1,\dots,\overline{k}_n)$. It follows that $l_{uw}=0$ if
$u\neq w$, and $l_{uu}=\frac{\lambda_j}{\lambda_i}$ for all
$u=1,\dots,n$. Hence the identity (\ref{cequ6}) is equivalent to
$\sum_{v=1}^md_vr_{d_v}a_u^{-1}=\frac{\lambda_j}{\lambda_i}r_{a_u}$.
\qed

\vspace{5mm}
\subsection*{Acknowledgement} {\small The work is
supported in part by an FWO-grant and NSFC (No. 10801099). The first
named author wishes to thank Department of Algebra and Analysis of
 University of Almeria for its hospitality during his visit.}

\begin{appendix}
\section{}

In this appendix, we list some properties of the cohom functors of
categories of comodules. These properties are probably well known.
Since we could not find any reference, we give a complete account of
proofs here.

Let $C$ be an arbitrary coalgebra. If ${}^C\!M$ is a quasi-finite
comodule, there is a cohom functor $h_C(M,-):{}^C\M\to Vect_\k$. The
cohom functor $h_C(M,-)$ is left adjoint to the tensor functor, that
is; for any left $C$-comodule $X$ and any vector space $V$, we have
a natural isomorphism (cf. \cite{Tak})
$$\Phi_{X,V}:\Hom(h_C(M,X),V)\longrightarrow\Hom_C(X,M\ot V).$$ Let $\eta$
be the unit of the adjoint pair $(h_C(M,-),M\ot-)$.

Let $D$ be another coalgebra. If $X$ is a $C$-$D$-bicomodule, then
$h_C(M,X)$ is a right $D$-comodule. Regard $h_C(M,X)$ as a left
$D^*$-module. Then $\Hom(h_C(M,X),V)$ is a right $D^*$-module.
Simultaneously, $\Hom_C(X,M\ot V)$ is also a right $D^*$-module with
the right $D^*$-action induced by the left $D^*$-action on $X$.

\begin{lem}\label{alem1} The natural isomorphism $\Phi_{X,V}$ is right $D^*$-module
isomorphism.
\end{lem}
\proof We have to show $\Phi_{X,V}(f\cdot d^*)=\Phi_{X,V}(f)\cdot
d^*$ for all $f\in\Hom(h_C(M,X),V)$ and $d^*\in D^*$. Note that
$\Phi_{X,V}(f)$ is the composition $\xymatrix@C=1cm{X
\ar[r]^{\eta_X\qquad\ }&M\ot h_C(M,X) \ar[r]^{\qquad id\ot f}&M\ot
V}$. We use $\rho^D_X$ to denote the right $D$-comodule structure
map of $X$, and use $\rho_{h_C(M,X)}^D$ to denote the right
$D$-comodule structure map of $h_C(M,X)$. We have the following
commutative diagram:
$$\xymatrix@C=1.5cm{
  X \ar[d]_{\rho_X^D} \ar[r]^{\eta_X\qquad} & M\ot h_C(M,X) \ar[d]_{id\ot\rho_{h_C(M,X)}^D}  &  \\
  X\ot D \ar[d]_{id\ot d^*} \ar[r]^{\eta_X\ot id\qquad\ } &M\ot h_C(M,X)\ot D  \ar[d]_{id\ot id\ot d^*}
  \ar[r]^{\ \qquad id\ot f\ot d^*} & M\ot V\ar[d]^{=} \\
 X \ar[r]^{\eta_X\qquad\ }&M\ot h_C(M,X)
\ar[r]^{\qquad id\ot f}&M\ot V.}$$ So, we obtain
$$\begin{array}{ccl}
(\Phi_{X,V}(f)\cdot d^*)(x)&=&\Phi(f)\circ(id\ot
d^*)\circ\rho_X^D(x)\\
&=&(id\ot f)\circ(\eta_X)\circ(id\ot d^*)\circ\rho_X^D(x)\\
&=&(id\ot f)\circ(id\ot id\ot
d^*)\circ(id\ot\rho_{h_C(M,X)}^D)\circ\eta_X(x)\\
&=&(id\ot f)\circ(id\ot f\cdot d^*)\circ\eta_X(x)\\
&=&\Phi_{X,V}(f\cdot d^*)(x).
 \end{array}
$$
Hence $\Phi_{X,V}$ is a right $D^*$-module morphism. \qed

\begin{lem}\label{alem2} Let ${}^CM$ and ${}^CN$ be quasi-finite
comodules, and $f:M\to N$ be a comodule morphism. Let ${}^CX^D$ be a
bicomodule. We have a commutative diagram of right $D^*$-module
morphisms:
$$\xymatrix@C=2.5cm{
  \Hom(h_C(M,X),V) \ar[d]_{\Phi_{M,V}} \ar[r]^{\Hom(h_C(f,X),V)} & \Hom(h_C(N,X),V) \ar[d]_{\Phi_{N,V}}  \\
  \Hom_C(X,M\ot V) \ar[r]^{\Hom_C(X,f\ot V)} & \Hom_C(X,N\ot V). }$$
\end{lem}
\proof The diagram follows from the following commutative diagrams in which the morphisms are natural ones: {\small$$\xymatrix@C=0.5cm{
 \Hom(\underset{\rightarrow_\lambda}\lim\Hom_C(X_\lambda,M)^*,V)
   \ar[d] \ar[r] & \underset{\leftarrow_\lambda}\lim\Hom(\Hom_C(X_\lambda,M)^*,V)
   \ar[d]\ar[r]&\underset{\leftarrow_\lambda}\lim\Hom_C(X_\lambda,M\ot V) \ar[d]\\
 \Hom(\underset{\rightarrow_\lambda}\lim\Hom_C(X_\lambda,N)^*,V)
 \ar[r] & \underset{\leftarrow_\lambda}\lim\Hom(\Hom_C(X_\lambda,N)^*,V)
 \ar[r]&\underset{\leftarrow_\lambda}\lim\Hom_C(X_\lambda,N\ot V)
 ,}$$}

$$\xymatrix@C=1cm{\Hom(h_C(M,X),V) \ar[d] \ar[r]^{\cong\qquad} &\Hom(\underset{\rightarrow_\lambda}\lim\Hom_C(X_\lambda,M)^*,V)
   \ar[d]\\
   \Hom(h_C(N,X),V) \ar[r]^{\cong\qquad} &\Hom(\underset{\rightarrow_\lambda}\lim\Hom_C(X_\lambda,N)^*,V),}$$
$$\xymatrix@C=1cm{\underset{\leftarrow_\lambda}{\lim}\Hom_C(X_\lambda,M\ot V) \ar[d]\ar[r]^{\ \cong}&\Hom_C(X,M\ot V)\ar[d]\\
\underset{\leftarrow_\lambda}\lim\Hom_C(X_\lambda,N\ot V) \ar[r]^{\
\cong}&\Hom_C(X,N\ot V).}$$
In the above diagrams, the limits run
through all the finite dimensional left $C^*$-subcomodule of $X$.
\qed

Let $e\in C^*$ be an idempotent. Let ${}^CY$ be a comodule. Then
$Ye$ is a left $eCe$-comodule. We have a left $C$-comodule morphism (cf. \cite{CG}):
\begin{equation}\label{atag2}
\theta_Y:Y\to eC\square_{eCe}Ye,\ y\mapsto \sum_{(y)}ey_{(-1)}\ot
y_{(0)}e,
\end{equation}
 where $\sum_{(y)}y_{(-1)}\ot y_{(0)}=\rho(y)$.
Since $eC$ is a quasi-finite left $C$-comodule, we have a natural
isomorphism
\begin{equation}\label{atag1}
\Psi:\Hom_C(Y,eC\square_{eCe}Ye)\longrightarrow\Hom_{eCe}(h_C(eC,Y),Ye).\end{equation}
Now let $\xi_Y=\Phi(\theta_Y):h_C(eC,Y)\longrightarrow Ye$.
If $Y$ is a $C$-$D$-bicomodule, then $\theta_Y$ is a
morphism of right $D$-comodules. One may check that
$\xi_Y=\Phi(\theta_Y)$ is also a morphism of right $D$-comodules.

\begin{lem}\label{alem3} The map $\xi_Y$ defined above is an
isomorphism. Moreover, we have a natural isomorphism
$$\xi:h_C(eC,-)\to (-)e$$ of functors from the category of left
$C$-comodules to the category of left $eCe$-comdules.
\end{lem}
\proof We already know from \cite{CG} that the functor $h_C(eC,-)$
is natural isomorphic to $(-)e$. We need to show that $\xi$ is
exactly the natural isomorphism between these two functors. From the
proof of \cite[Theorem 1.5]{CG}, we know that there is a natural
isomorphism $$\Psi:\Hom_C(Y,eC\square_{eCe}Z)\longrightarrow
\Hom_{eCe}(Ye,Z),$$ for $Y\in {}^C\M$ and $Z\in {}^{eCe}\M$.
Moreover, given $f\in\Hom_C(Y,eC\square_{eCe}Z)$, we have
$\Psi(f)=(e\ot id)\circ f$. Now let $Z=Ye$. One sees that
$\Psi(\theta_Y)=id_{Ye}$. So, $\theta_Y:Y\to eC\square_{eCe}Ye$ is
the unit map of the adjoint pair $((-)e,eC\square_{eCe}-)$. Since
$(h_C(eC,-),eC\square_{eCe}-)$ is also an adjoint pair, the unit map
$\theta_Y$ induces an isomorphism through the isomorphism
(\ref{atag1}). That is, $\xi_Y:h_C(eC,Y)\longrightarrow Ye$ is a
natural isomorphism. \qed

\begin{prop}\label{aprop} Let $X$ be a left $C$-comodule, and $e_1,e_2\in C^*$ be
idempotents. Given an element $c^*\in e_2C^*e_1$, we have a left
$C$-comodule morphism $e_1C\overset{c^*\cdot}\longrightarrow e_2C$
and a commutative diagram:
$$\xymatrix@C=2cm{h_C(e_2C,X)\ar[d]_{\xi^2_X}\ar[r]^{h_C(c^*\cdot,X)}&h_C(e_1C,X)\ar[d]^{\xi^1_X}\\
Xe_2\ar[r]^{\cdot c^*}&Xe_1,}$$ where $\xi^1_X$ and $\xi^2_X$ are
natural isomorphisms formed in Lemma \ref{alem3} corresponding to
idempotents $e_1$ and $e_2$ respectively.
\end{prop}
\proof To show the diagram to be commutative, it suffices to
prove that the following compositions of morphisms coincide:
$$\xymatrix{\varphi:M\ar[r]^{\eta^2_M\qquad\quad}&e_2C\ot h_C(e_2C,X)\ar[r]^{\qquad id\ot\xi_X^2}& e_2C\ot Xe_2\ar[r]^{id\ot \cdot c^*}& e_2C\ot Xe_1,}$$
$$\xymatrix{\phi:M\ar[r]^{\eta^2_M\qquad\quad}&e_2C\ot h_C(e_2C,X)\ar[rr]^{\qquad id\ot h_C(c^*\cdot,X)\qquad}
&& e_2C\ot h_C(e_1C,X)\ar[r]^{\ \quad id\ot \xi_X^1}& e_2C\ot
Xe_1,}$$ where $\eta_M^2$ is the unit of the corresponding adjoint
pair. By Lemma \ref{alem3}, we have $\varphi=(id\ot \cdot
c^*)\circ\theta^2_X$, where $\theta^2_X$ is formed in (\ref{atag2})
corresponding to the idempotent $e_2$. In the commutative diagrams of
Lemma \ref{alem2}, if we set $N=e_1C$, $M=e_2C$ and $V=h_C(e_2C,X)$,
then we obtain the following commutative diagram
$$\xymatrix@C=2cm{X\ar[d]_{\eta_M^2}\ar[r]^{\eta^1_M\qquad}&e_1C\ot h_C(e_1C,X)\ar[d]^{c^*\cdot\ot id}\\
e_2C\ot h_C(e_2C,X)\ar[r]^{id\ot h_C(c^*\cdot, X)}&e_2C\ot
h_C(e_1C,X).}$$ Hence
$$\begin{array}{ccl}
\phi&=&(id\ot \xi_X^1)\circ(c^*\cdot\ot
id)\circ\eta_X^1\\
&=&(c^*\cdot\ot id)\circ(id\ot \xi_X^1)\circ\eta_X^1\\
&=&(c^*\cdot\ot id)\circ\theta_X^1,
  \end{array}
$$ where $\theta_X^1$ is formed in (\ref{atag2}) corresponding to
the idempotent $e_1$. Now for any element $x\in X$, we have
$$\begin{array}{ccl}\varphi(x)&=&\sum_{(x)}e_2x_{(-1)}\ot
x_{(0)}e_2c^*\\
&=&\sum_{(x)}e_2x_{(-1)}\ot
x_{(0)}e_2c^*e_1\\
&=&\sum_{(x)}x_{(-5)}e_2(x_{(-4)})\ot
e_2(x_{(-3)})c^*(x_{(-2)})e_1(x_{(-1)})x_{(0)}\\
&=&\sum_{(x)}x_{(-4)}e_2(x_{(-3)})\ot
c^*(x_{(-2)})e_1(x_{(-1)})x_{(0)},  \end{array}$$ and
$$\begin{array}{ccl}\phi(x)&=&\sum_{(x)}c^*e_1x_{(-1)}\ot
x_{(0)}e_1\\
&=&\sum_{(x)}e_2c^*e_1x_{(-1)}\ot
x_{(0)}e_1\\
&=&\sum_{(x)}x_{(-4)}e_2(x_{(-3)})\ot
c^*(x_{(-2)})e_1(x_{(-1)})x_{(0)}.
\end{array}$$ Hence $\varphi=\phi$ as desired. \qed

\begin{rem} If $X$ is a $C$-$D$-bicomodule, then the morphisms in the commutative
diagram of last proposition are morphisms of right $D$-comodules.
\end{rem}

\end{appendix}

\bibliography{}

\end{document}